\documentclass[12pt]{article}
\usepackage{amsmath,amsthm,amsfonts,amssymb,amscd,enumerate}
\RequirePackage{graphics}

 \theoremstyle{plain}

\newtheorem*{mainthm}{The Main Theorem}
 
 \theoremstyle{definition}
 \theoremstyle{remark}

\newtheorem*{remarks}{Remarks}

\newcommand{\zz}{{\mathbb Z}}

\newcommand{\minus}{{-1}}

 \newcommand{\gb}{\beta}

\newcommand{\gs}{\sigma}

\newcommand{\theh}{\theta_3^H}

 \newcommand{\gD}{\Delta}

\newcommand{\Ker}{\operatorname{Ker}}

\newcommand{\sq}{\operatorname{Sq^1}}

\newcommand{\cat}{{\mathrm{CAT}}}
\newcommand{\res}{{\mathrm{res}}}

\newcommand\mapright[1]{\,\smash{\mathop{\longrightarrow\,}\limits^{#1}}}

    {
\end{enumerate}
}
\newenvironment{enumerate1*}{
\begin{enumerate}[\upshape (*1)]}%
	{
\end{enumerate}
}
\newenvironment{enumeratea}{
\begin{enumerate}[\upshape (a)]}%
	{
\end{enumerate}
}
\newcommand{\comment}[1]{}

\usepackage{color}
\usepackage[outerbars,color]{changebar}
\usepackage{ifthen}

\newboolean{usecolor}
\setboolean{usecolor}{true}

\ifthenelse{\boolean{usecolor}}
{
  \definecolor{colore}{cmyk}{0,1,0.6,0}
  \definecolor{coloregen}{cmyk}{0.7,0,1,0}
  \definecolor{coloresimo}{cmyk}{1,0.6,0,0}
}
{
  \definecolor{colore}{cmyk}{0,0,0,1}
  \definecolor{coloregen}{cmyk}{0,0,0,1}
  \definecolor{coloresimo}{cmyk}{0,0,0,1}
}

		{\end{enumerate}}

\begin{document}

\title{Aspherical manifolds that cannot be triangulated}

\author{Michael W. Davis \thanks{The first author was partially supported by
NSF grant DMS 1007068.}    
\and Jim Fowler
\and Jean-Fran\c{c}ois Lafont \thanks{The third author was partially supported by NSF grant DMS 1207782. \newline 
Any opinions, findings, and conclusions or recommendations expressed in this material are those of the author(s) and do not necessarily reflect the views of the National Science Foundation.}}
\date{\today} \maketitle

\begin{abstract}
Although Kirby and Siebenmann showed that there are manifolds which do not admit  PL structures, the possibility 
remained that all manifolds could be triangulated.  
In the late seventies Galewski and Stern constructed a closed $5$-manifold $M^5$ so that every $n$-manifold, with $n\ge 5$, can be triangulated if and only if $M^5$ can be triangulated.  Moreover, $M^5$ admits a triangulation if and only if the Rokhlin $\mu$-invariant homomorphism, $\mu:\theh \to \zz/2$, is split. In 2013 Manolescu showed that the $\mu$-homomorphism does \emph{not} split. Consequently, there exist  Galewski-Stern manifolds, $M^n$, that are not triangulable for each $n\ge 5$.  In 1982 Freedman proved that there exists a topological $4$-manifold with even intersection form of signature 8.  It  followed from later work of Casson that such $4$-manifolds cannot be triangulated.  In 1991 Davis and Januszkiewicz applied Gromov's hyperbolization procedure to Freedman's $E_8$-manifold to show that there exist closed aspherical $4$-manifolds that cannot be triangulated.  In this paper we apply hyperbolization techniques to the Galewski-Stern manifolds  to show that there exist closed  aspherical $n$-manifolds that cannot be triangulated for each $n\ge 6$.  The question remains open in dimension $5$.
\smallskip

\noindent
\textbf{AMS classification numbers}.  Primary: 57Q15, 57Q25, 20F65,\\
Secondary: 57R58
\smallskip

\end{abstract}

\section{Introduction}\label{intro}
Any $3$-dimensional homology sphere, $H^3$, bounds a PL $4$-manifold, $W^4$, with vanishing first and second Stiefel-Whitney classes.  The intersection form of $W^4$ is then unimodular and even; so, its signature, $\gs(W^4)$, is divisible by $8$.  If $H^3$ is the standard $3$-sphere,  then, by Rokhlin's Theorem, $\gs(W^4)$ is divisible by $16$.  So, one defines the $\mu$-invariant of $H^3$ by
\[
\mu(H^3)= \frac{1}{8} \gs(W^4) \mod 2.
\]
Let $\theh$ be the abelian group (under connected sum) of homology cobordism classes oriented PL homology $3$-spheres.  The $\mu$-invariant  defines a homomorphism $\mu:\theh \to \zz/2$. (For background about this material, see \cite{sav}.)

In \cite{ks} Kirby and Siebenmann proved that for any topological manifold $M^n$ there is an obstruction $\gD\in H^4(M^n;\zz/2)$ which, for $n\ge 5$, vanishes if and only if $M^n$ admits a PL structure.  An important  point here is that the Kirby-Siebenmann obstruction can be defined for any polyhedral homology manifold $M^n$, as follows.  First, there is obstruction in $H^4(M^n;\theh)$ to finding an acyclic resolution of $M$ by a PL manifold.  This is the class of the cocycle that associates to each $4$-dimensional ``dual cell'' in $M^n$ the class of its boundary in $\theh$.  The Kirby-Siebenmann obstruction $\gD$ is the image of this element of $H^4(M^n;\theh)$ under the coefficient homomorphism $\mu:\theh \to \zz/2$.

After the proof by Edwards and Cannon of the Double Suspension Theorem,  it seemed possible that every  topological manifold could still be homeomorphic to some simplicial complex, even if it did not have a PL structure. Galewski-Stern \cite{gsannals, gsams} and independently Matumoto \cite{matumoto} proved that the following statements are equivalent:
\begin{enumeratea}
\item
Every manifold of dimension $\ge 5$ can be triangulated.
\item
There exists a homology $3$-sphere $H^3$ with $\mu(H^3)=1$ and $[H^3]$ of order $2$ in $\theh$.
\end{enumeratea}
Galewski-Stern also showed that, for $n\ge 5$, the obstruction for a manifold to have a simplicial triangulation was the Bockstein of its Kirby-Siebenmann obstruction, $\gb(\gD)\in H^5(M^n;\Ker \mu)$, where $\gb:H^4(M^n;\zz/2) \to H^5(M^n;\Ker \mu)$ denotes the Bockstein homomorphism associated to the short exact sequence of coefficients,
\begin{equation}\label{e:bockstein}
    0\mapright{} \Ker \mu \mapright {} \theh \mapright {\mu} \zz/2 \mapright{} 0 \ .
	\end{equation}
In fact, as was shown in \cite{gsuniversal}, one can focus instead on a simpler Bockstein associated to the short exact sequence ofcoefficients, 
\begin{equation}\label{e:usual}
    0\mapright{} \zz/2\mapright{\times 2}\, \zz/4 \mapright{} \zz/2\mapright{} 0\ .
    \end{equation}
As is well-known, the Bockstein associated to \eqref{e:usual} is the first Steenrod square, $\sq$.  So, henceforth we will use $\sq$ instead of a more general Bockstein $\gb$.  The reduction to the case of $\sq$ goes as follows.  Suppose $N^n$ is a manifold with $\sq(\gD)\neq 0$.  By \cite[Theorem 2.1]{gsuniversal}, if one such $N^n$ can be triangulated, then every manifold of dimension $\ge 5$ can be triangulated (this is statement (a) above).  In \cite{gsuniversal} Galewski-Stern also constructed $n$-manifolds,  for each $n\ge 5$, with $\sq(\gD)\neq 0$.

Manolescu has recently established that homology $3$-spheres satisfying (b) do not exist \cite[Cor.~1.2]{mano}.  It follows that any manifold with $\sq(\gD) \ne 0$ is not homeomorphic to a simplicial complex.  So the Galewski-Stern manifolds cannot be triangulated.

By work of Freedman and Casson nontriangulable manifolds exist in dimension 4 (cf.~\cite{am}).  First, Freedman \cite{freedman} showed that any homology $3$-sphere bounds a contractible (topological) $4$-manifold.  One defines the $E_8$-homology manifold $X^4$ as follows. Start with the plumbing $Q(E_8)$ defined by the $E_8$ diagram.  It is a smooth, parallelizable $4$-manifold with boundary; its boundary being Poincar\'e's homology $3$-sphere, $H^3$. The signature of $Q(E_8)$ is $8$.  $X^4$ is defined to be the union of this plumbing with $c(H^3)$ (the cone on $H^3$).  It is a polyhedral homology $4$-manifold with one non-manifold point.  By \cite{freedman} we can topologically ``resolve the singularity'' of $X^4$ by replacing $c(H^3)$ with a contractible $4$-manifold bounded by $H^3$ to obtain a topological $4$-manifold $M^4$ with nontrivial Kirby-Siebenmann invariant.  $M^4$ cannot be triangulated.  (Proof: For any  triangulation of $M^4$, the link of a vertex is a simply connected $3$-dimensional homology sphere; so,  by Perelman's proof of the Poincar\'e Conjecture, it is $S^3$. Alternatively, since the Casson invariant of any homotopy $3$-sphere is $0$, so is its $\mu$-invariant.   So any triangulation of $M^4$ would automatically be PL, contradicting Rokhlin's Theorem.) A variation of this idea was used 
in \cite{dj} to produce an aspherical $4$-manifold that cannot be triangulated.  Start with a triangulation of $X^4$.  Apply the technique of \cite {dj} to $X^4$ to get its ``hyperbolization''  $h(X^4)$.  It is a $\cat(0)$ polyhedral homology manifold with one non-manifold point, the link of which is $H^3$.  Resolve the conical singularity to obtain a closed aspherical topological manifold $N^4$ with $\gD(N^4)\neq 0$. By the previous argument, $N^4$ cannot be triangulated. (Of course, $N^4\times S^1$ can be triangulated since it is homeomorphic to the triangulated manifold $X^4\times S^1$.)

The idea of this paper is to apply hyperbolization techniques to the Galewski-Stern manifolds to obtain aspherical manifolds $N^n$ that cannot be triangulated.  We do not know how to make our techniques work in dimension 5; however, they do work in any dimension $\ge 6$.  So, we get the following.  
\begin{mainthm}
For each $n\ge 6$ there is a closed aspherical manifold $N^n$ that cannot be triangulated.
\end{mainthm}

Our thanks go to Ron Stern for some helpful comments. 

\section{The construction}
\paragraph{The Galewski-Stern $\boldsymbol{5}$-manifold.}
We recall the Galewski-Stern construction in \cite{gsuniversal} of a $5$-manifold, $N^5$, now known to be nontriangulable.  Start with $X^4 \times I$, where $X^4$ is the $E_8$-homology manifold. Attach an orientation-reversing $1$-handle, $D^3\times I$, connecting the two copies of $c(H^3)$ along their boundary.  The two copies of $Q(E_8)$ become the boundary connected sum $Q(E_8)\#_b\, Q(E_8)$; it is  a  $4$-manifold with boundary, the boundary being $H^3 \# H^3$  (not $H^3\# (-H^3)$).  Consider $c(H^3) \cup  (D^3\times I) \cup\, c(H^3)$.  It is a contractible polyhedral homology $4$-manifold; its boundary is $H^3\#H^3$.  Fill in the boundary with $c(H^3\# H^3)$ to obtain a polyhedral homology manifold $T$ with the homology of $S^4$ (i.e., a ``generalized homology $4$-sphere'').  Next fill in $T$ with $c(T)$.  The result is the polyhedral $5$-manifold  with boundary
    \begin{equation}\label{e:p}
    P^5:=(X^4 \times I)  \cup c(T).
    \end{equation}
Roughly speaking, after ignoring the differences between homology spheres and spheres, we have attached an orientation-reversing $1$-handle $D^4\times I$ to $X^4\times I$ (the $1$-handle is actually $c(T)$). So, $P^5$ is a nonorientable, polyhedral homology $5$-manifold with boundary; it is homotopy equivalent to $S^1\vee X^4$.  The boundary of $P^5$ is $c(H^3\#H^3) \cup \,Q(E_8)\#_b\, Q(E_8)$; so, $\partial P^5$ contains a single non-manifold point (the cone point of $c(H^3\#H^3)$).
By Edwards' Polyhedral-Topological Manifold Characterization Theorem  \cite[p.~119] {edwards}, the interior of $P^5$ is a topological manifold.  Its Kirby-Siebenmann invariant, $\gD(P^5)$, is the image of $[X^4]$ in $H^4(P^5;\zz/2)$. Thus, $P^5$ is polyhedral homology manifold with boundary and its interior is a topological manifold which does not admit a PL structure.    Since the Kirby-Siebenmann obstruction of $\partial P$ is 0, $\gD(P)$ is the image of a (unique) class $\gD(P,\partial P) \in H^4(P,\partial P;\zz/2)$; moreover,  one sees  that the image $\sq(\gD(P,\partial P))$ under $\sq$ is the nonzero class in $H^5(P,\partial P;\zz/2)$. 

When it comes to  applying hyperbolization, it is at this point where the Galewski-Stern construction becomes problematic.  Galewski and Stern get rid of the singular point of $\partial P^5$ as follows:  (1) attach an external collar $\partial P \times [0,1]$ to $P$, (2) find a PL manifold $V^4$ embedded in $\partial P \times (0,1)$, (3) define $U$ to be the part of the external collar between $\partial P \times 0$ and $V^4$, (4) argue that $V^4$ bounds a PL 5-manifold $W$ (necessarily nonorientable), and finally, (5) glue in $W$ to get the desired manifold, 
\[
N^5:= P \cup U\cup W .
\]

\paragraph{Galewski-Stern manifolds of dimension $\boldsymbol{> 5}$.}
In dimensions $>5$ there are versions of these manifolds to which it is easier to apply hyperbolization.  To fix ideas, suppose the dimension is $6$.  Let $P':= P^5\times S^1$, where $P^5$ is the polyhedral homology $5$-manifold with boundary which was constructed previously. Then $\gD(P')$ is the nontrivial element in $H^4(P;\zz/2)\otimes H^0(S^1;\zz/2)$, a summand of $H^4(P';\zz/2)$.  It is the image of the unique nontrivial element $\gD(P',\partial P')$ of $H^4(P,\partial P;\zz/2)\otimes H^0(S^1;\zz/2)$. Its image under $\sq$ is denoted $\sq(\gD(P',\partial P')) \in H^5(P,\partial P;\zz/2)\otimes H^0(S^1;\zz/2)\cong \zz/2$.  By Edwards' Theorem, $\partial P'$ is a topological $5$-manifold. Since $\gD(\partial P')$ is zero, $\partial P'$ is actually homeomorphic to a PL $5$-manifold $V'$. It is easy to see that $V'$ bounds a PL $6$-manifold $W'$.  Put
\[
N':= P' \cup U \cup W'
\]
where $U$ is the mapping cylinder of a homeomorphism between $\partial P'$ and
$V'=\partial W'$.
Since $\gD(N')$ restricts to $\gD(P')$, we have $\gD(N')\neq 0$ and 
one sees as before that $\gD(N')$ is the image of 
\[
\gD(P',\partial P')\in H^4(P',\partial P';\zz/2)\cong H^4(N',U\cup W';\zz/2).
\]
By Wu's Formula, $\sq(\gD(P',\partial P'))=w_1(P')\cup \gD(P',\partial P')\neq 0$. (This argument is from the final paragraph of \cite{gsuniversal}.)  Hence, $\sq(\gD(N'))$ is the nonzero image of $\sq(\gD(P',\partial P'))$.

\paragraph{Hyperbolization.} A hyperbolization technique of Gromov \cite{gromov} is explained in \cite{dj}:  given a simplicial complex $K$, one can construct a new space $h(K)$ and a map $f:h(K)\to K$ with the following properties.
\begin{enumeratea}
\item
$h(K)$ is a locally $\cat(0)$ cubical complex; in particular, it is aspherical.
\item
The inverse image in $h(K)$ of any simplex of $K$ is a ``hyperbolized simplex''.  So, the inverse image of each vertex in $K$ is a point in $h(K)$.
\item
$f:h(K)\to K$ induces a split injection on cohomology (cf.~\cite[p.~355]{dj}).
\item 
Hyperbolization preserves local structure: for any simplex $\gs$ in $K$ the link of $f^\minus(\gs)$ is isomorphic to a subdivision of the link of $\gs$ in $K$ (cf.~\cite[p.~356]{dj}).  So, if $K$ is a polyhedral homology manifold, then so is $h(K)$.
\item
If $K$ is a polyhedral homology manifold, then $f:h(K)\to K$ pulls back the Stiefel-Whitney classes of $K$ to those of $h(K)$. 
\end{enumeratea}

In \cite{djw} the above version of hyperbolization is used to define a ``relative hyperbolization procedure'' (an idea also due to Gromov \cite{gromov}).  Given $(K,\partial K)$, a triangulated manifold with boundary, form $K\cup c(\partial K)$ and then define $H(K,\partial K)$ to be the complement of an open neighborhood of the cone point in $h(K\,\cup\, c(\partial K))$.  Then $H(K,\partial K)$ is a manifold with boundary; its boundary is homeomorphic to $\partial K$.  The main result of \cite{djw} is that if each component of $\partial K$ is aspherical, then so is $H(K,\partial K)$; moreover,  the inclusion $\partial K\to H(K,\partial K)$ induces an injection on fundamental groups for each component of $\partial K$.  In other words, if a triangulated aspherical manifold bounds a triangulated manifold, then it bounds an aspherical manifold.

\begin{proof}[Proof of the Main Theorem] Our nontriangulable $6$-manifold $N^6$ will be the union of three pieces, $N^6=P_1\cup U\cup P_2$, where $P_1$ and $P_2$ are  triangulable, aspherical $6$-manifolds with boundary, and where $U$ is the mapping cylinder of a homeomorphism $\partial P_1 \to \partial P_2$.  $P_1$ will be defined via hyperbolization and $P_2$ via relative hyperbolization.  Put
\(
P_1:=h(P'),
\)
where $P'=P^5\times S^1$ is the simplicially triangulated $6$-manifold with boundary discussed above.  Then $P_1$ is a topological $6$-manifold with boundary and $\partial P_1 =h(\partial P')$ is homeomorphic to a PL $5$-manifold $V$.  (N.B. the PL structure on $V$ is incompatible with the  triangulation of $\partial P_1$ as a subcomplex of $P_1$.) Let $U$ be the mapping cylinder of a homeomorphism $h(\partial P')\to V$.  Let $W$ be a PL manifold bounded by $V$.  Equip $W$ with a PL triangulation.  Applying relative hyperbolization, we get an aspherical $6$-manifold with boundary $P_2:=H(W,V)$; its boundary being $V$.  Then  $N^6=P_1\cup U\cup P_2$ is aspherical.  
By properties (c) and (d) of hyperbolization, $\gD(P_1,\partial P_1)=f^*(\gD(P',\partial P'))$. So, $\sq(\gD(P_1,\partial P_1))=f^*(\sq(\gD(P',\partial P')))$.  Consequently,  $\gD(P_1,\partial P_1)$ and $\sq(\gD(P_1,\partial P_1))$ are both nonzero.  Since $P_2$ is a PL manifold, its obstructions vanish. As before, it follows that $\gD(N^6)$ and $\sq(\gD(N^6))$ are both nonzero.
\end{proof}

\begin{remarks}
(a) What is the situation in dimension $5$? As explained in the Introduction,  any polyhedral homology $4$-manifold, $M^4$, can be resolved to a topological manifold $M^4_\res$.   
When $P^5$ is defined by \eqref{e:p}, $(\partial P)_\res$ does not support a PL structure (although $\gD((\partial P)_\res)=0$). However, one can vary the definition of $P$ so that $(\partial P)_\res$ becomes homeomorphic to a PL $4$-manifold.  To complete our construction we need $(h(\partial P))_\res$, the resolution of the hyperbolization, to be PL.  If this could be achieved, we could finish by finding a PL $5$-manifold bounded by $(h(\partial P))_\res$ and applying relative hyperbolization as before.
 
(b) The Galewsi-Stern $5$-manifold is nonorientable.  (In fact, Siebenmann showed in \cite{sieb} that every orientable $5$-manifold can be triangulated.)  Since our construction of a nontriangulable aspherical manifold $N^6$ starts from $P^5\times S^1$,  the manifold $N^6$ is also nonorientable (by property (e) of hyperbolization).  The question arises: do orientable examples exist?  The answer is yes. To construct one, start from the nonorientable $S^1$-bundle over $P$ instead of $P\times S^1$, where $w_1$ of the associated vector bundle is $w_1(P)$.  The total space $E$ of the $S^1$-bundle is then a $6$-dimensional, orientable, polyhedral homology manifold with boundary. The restriction of the $S^1$-bundle to $\partial P$ is orientable so $\partial E=\partial P\times S^1$ is the same. As before, we get $N^6=P_1\cup U\cup P_2$ where $P_1=h(E)$, $V$ is a PL-manifold manifold homeomorphic to $\partial P_1$, $W$ is a PL manifold bounded by $V$, $P_2=H(W,V)$, and $U$ is the mapping cylinder of a homeomorphism  $\partial P_2\to \partial P_1$. 

(c) By using relative hyperbolization, it is proved in \cite{djw} that if an aspherical manifold bounds a triangulable manifold, then it bounds an aspherical one.  Can we omit the word ``triangulable?''
In other words,
if an aspherical topological manifold M bounds, does it bound an aspherical
manifold?    
If M does not support any triangulation, then it
cannot bound a triangulated manifold, and one can no longer 
use relative hyperbolization directly. For the specific examples
of nontriangulable aspherical manifolds given in this
paper, a similar construction can be applied to produce aspherical manifolds that they bound. The general case remains open.
\end{remarks}

\obeylines
Michael W. Davis, Department of Mathematics, The Ohio State University, 231 W. 18th Ave., Columbus Ohio 43210  {\tt davis.12@math.osu.edu} 
Jim Fowler, Department of Mathematics, The Ohio State University, 231 W. 18th Ave., Columbus Ohio 43210  {\tt fowler@math.ohio-state.edu} 
Jean-Fran\c{c}ois Lafont, Department of Mathematics, The Ohio State University, 231 W. 18th Ave., Columbus Ohio 43210  {\tt jlafont@math.ohio-state.edu}

\end{document}